\documentclass[11pt]{article}
\usepackage{lmodern} 
\usepackage[T1]{fontenc}
\usepackage{textcomp}
\usepackage[scr=boondoxo,frak=boondox,bb=boondox]{mathalfa}
\usepackage{epsfig}
\usepackage{t1enc}
    \usepackage[latin1]{inputenc}
    \usepackage[english]{babel}
        \usepackage{latexsym}
\usepackage{amssymb}
\usepackage{euscript}
\def\ep{\varepsilon}
\begin{document}
\setlength{\arraycolsep}{.136889em}
 \renewcommand{\theequation}{\thesection.\arabic{equation}}
\newcommand{\dto}{\stackrel{{d}}{\rightarrow}}
 \newtheorem{thm}{Theorem}[section]
 \newtheorem{propo}{Proposition}[section]
 \newtheorem{lemma}{Lemma}[section]
 \newtheorem{corollary}{Corollary}[section]
 \newtheorem{remark}{Remark}[section]
 \def\begg{\begin{equation}}
 \def\endd{\end{equation}}
 \def\ep{\varepsilon}
 \def\noo{n\to\infty}
 \def\al{\alpha}  
 \def\be{\bf E}
 \def\bp{\bf P}
 \medskip
 \centerline{\Large\bf On the local time of anisotropic random walk on $\mathbb Z^2$}

 \bigskip\bigskip

 \centerline{\large Dedicated to the memory of P\'al R\'ev\'esz} 
 
  \smallskip
 
  \centerline{\large  who was our inspiration, collaborator and friend for half a century.}
  
   \bigskip\bigskip
   \noindent
 {\bf Endre Cs\'aki}

 \noindent
 Alfr\'ed R\'enyi Institute of Mathematics, Budapest, P.O.B. 127, H-1364, Hungary. E-mail address:
 csaki.endre@renyi.hu

 \bigskip
 
 \bigskip
  \noindent
 {\bf Ant\'onia F\"oldes*}
 
  \noindent
 Department of Mathematics, College of Staten Island, CUNY, 2800
 Victory Blvd., Staten Island, New York 10314, U.S.A.  E-mail
 address: antonia.foldes@csi.cuny.edu*
 
 \bigskip\bigskip
 
  \noindent{\bf Abstract}\newline We study the local time of the
 anisotropic random walk on the two-dimensional lattice $\mathbb{Z}^2$, by establishing   the exact asymptotic behavior of the $N$- step  return probability to the origin.

 \medskip
 \noindent {\it MSC:} primary 60F17, 60G50, 60J65; secondary 60F15, 
 60J10.

 \medskip
 \noindent {\it Keywords:} anisotropic random walk;
strong approximation; 2-dimensional Wiener process; local time

 \section{Introduction and main results}
 \renewcommand{\thesection}{\arabic{section}} \setcounter{equation}{0}
 \setcounter{thm}{0} \setcounter{lemma}{0}

We consider random walks on the square lattice $\mathbb{Z}^2$ of the plane
with possibly unequal symmetric horizontal and vertical step
probabilities, so that these probabilities can only depend on the value of
the vertical coordinate. In particular, if such a random walk is situated
at a site on the horizontal line $y=j\in \mathbb{Z}$, then at the next
step it moves with probability $p_j$ to either vertical neighbor, and with
probability $1/2-p_j$ to either horizontal neighbor. A substantial
motivation for studying such two-dimensional random walks on anisotropic
lattice has originated from transport problems of statistical physics.

More formally, consider the random walk
$\{{ \mathbf C}(N)=\left(C_1(N),C_2(N)\right);\,\, N=0,1,2,\ldots\}$
on $\mathbb{Z}^2$ with the transition probabilities
$$ {\bf P}({\bf C}(N+1)=(k+1,j)|{\bf C}(N)=(k,j))={\bf P}
({\bf C}(N+1)=(k-1,j)|{\bf C}(N)=(k,j))=\frac{1}{2}-p_j,$$
$$ { \mathbf P} ({ \mathbf C}(N+1)=(k,j+1)|{ \mathbf C}(N)=(k,j))={ \mathbf P}
({ \mathbf C}(N+1)=(k,j-1)|{ \mathbf C}(N)=(k,j))=p_j,$$
for $(k,j)\in{\mathbb Z^2}$, $N=0,1,2,\ldots$ We assume throughout the
paper that $0<p_j\leq 1/2$ and $\min_{j\in\mathbb Z} p_j<1/2$.  We assume also that ${\mathbf C}(0)=(0,0)$.

In our paper \cite{CCFR13} we gave a detailed account of the history of this topic. Here we just mention some special cases and the corresponding references. 
The case $p_j=1/4,\, j=0,\pm 1,\pm 2,\ldots$ corresponds to the simple
symmetric random walk on the plane. For this case, we refer to Erd\H os and
Taylor \cite{ET}, Dvoretzky and Erd\H os \cite{DE}, and R\'ev\'esz \cite{RE}.
The case $p_j=1/2$ for some $j$ means that the horizontal line $y=j$ is
missing. If all $p_j=1/2$, then the random walk takes place on the $y$ axis,
so it is only a one-dimensional random walk, and this case is excluded from
the present investigations. The case however when $p_j=1/2,\,
j=\pm 1,\pm 2,\ldots$ but $p_0=1/4$ is an interesting one, which is the
so-called random walk on the two-dimensional comb. For the general model we may
refer to Weiss and Havlin \cite{WH}, Bertacchi and Zucca \cite{BZ},
Bertacchi \cite{BE}, Cs\'aki {\it et al.} \cite{CCFR08}. 
 For early investigations of the model we refer to
Silver {\it et al.} \cite{SSL}, Seshadri {\it et al.} \cite{SLS},  Shuler
\cite{SH}, Westcott \cite{WE}, where certain properties of this random
walk were studied under various conditions. Heyde \cite{H} proved an
almost sure approximation for $C_2(\cdot)$ under the  following condition 

\begg n^{-1}\sum_{j=1}^n p_j^{-1}=2\gamma+o(n^{-\eta}), \qquad
n^{-1}\sum_{j=1}^n p_{-j}^{-1}=2\gamma+o(n^{-\eta}),  \quad  {\rm as} \quad n\to \infty\label{h}   \endd 
for some constants $ \gamma< \infty$ and $1/2<\eta.$ 

He also proved that $C_2(N)$ is asymptotically normally distributed, namely  

\bigskip

\noindent
{\bf Theorem A} \cite{H} {\it Under {\rm(\ref{h})} with  $1/2<\gamma< \infty$ and} $1/2<\eta$
$$\gamma^{1/2}N^{-1/2}C_2(N)~\dto N(0,1), \quad {\rm as} \quad N\to\infty,$$
where $ ~\dto$ denotes convergence in distribution.

In \cite{CCFR13} we proved the following joint strong approximation result.

\bigskip
\noindent
{\bf Theorem B} \cite{CCFR13} {\it Under the condition {\rm (\ref{h})} with $1/2<\eta\leq 1$,
on an appropriate probability space for the random walk
$$\{{ \mathbf C}(N)=(C_1(N),C_2(N)); N=0,1,2,\ldots\}$$
one can construct two independent
standard Wiener processes $\{W_1(t);\, t\geq 0\}$, $\{W_2(t);\,
t\geq 0\}$ so that, as $N\to\infty$, we have with any}
$\varepsilon>0$
\begg
\left|C_1(N)-W_1\left(\frac{\gamma-1}{\gamma}\,N\right)\right|+
\left|C_2(N)-W_2\left(\frac{1}{\gamma}\,N\right)\right|
=O(N^{5/8-\eta/4+\varepsilon})\quad {a.s.}
\endd

Heyde {\it et al.} \cite{H} treated the case, when conditions
similar to (\ref{h}) are assumed but $\gamma$ can be different for the two
parts of (\ref{h}) and obtained almost sure convergence to
the so-called oscillating Brownian motion. The interested reader can find some other generalizations in Heyde \cite{H93}  and den Hollander \cite {DH}. Roerdink and Shuler \cite{RS} proved some asymptotic properties, including local limit theorems, under certain conditions. For more detailed history see \cite{DH} and   \cite{CF22B}.

In this paper we are interested in the local time of the anisotropic walk.
So far  we have investigated the following  particular cases. The local time of the comb was discussed  in our paper \cite{CSCSFR11}.  Another special case is the  so-called half-plane half-comb  walk;  when $p_j=1/4$, $j=0,1,2,\ldots$ and $p_j=1/2$, $j=-1,-2,\ldots$, i.e., we have a square lattice on the upper
half-plane and a comb structure on the lower half-plane. Its local time was given in \cite{CF22A}.
The anisotropic walk is called periodic with period $L$
 when $p_j=p_{j+L}$ for each $j\in \mathbb{Z},$
where $L\geq 1$ is a positive integer. The local time of this walk was investigated in \cite{CCFR13}. In that paper, the method of our proof was based on a result of Roerdink and Shuler \cite{RS}. Our aim in this paper is   to generalize this result for a much broader class of anisotropic walks, by using a totally different method. To get the local time one  usually finds first the limit of the $n$-step transition probability, and then gets local time results using the Green function.
Our main result is the following theorem.

\begin{thm} Consider the anisotropic walk defined above. Suppose that
\begin{itemize}
\item  {\rm(i)} \,\, {\rm(\ref{h})} holds with some    $1<\gamma< \infty$ and $1/2<\eta<\infty .$ 

\item {\rm(ii)}\,\, there is an $\omega>0 $ such that $p_j\geq \omega,\quad 
j=0,\pm 1,\pm 2,\ldots$ 

  \item {\rm(iii)} \,\, $ \sup_{n\geq 0}\sqrt{n} \sup_{N\geq n} \sum_{m=n}^N \left({ \mathbf P}(C_2(2m+2)=0 )-{ \mathbf P}(C_2(2m+1)=0 )\right)< \infty.$

\end{itemize}
Then 
$${ \mathbf P}({ \mathbf C}(2N)=(0,0))\sim \frac{1}{4Np_0\pi\sqrt{\gamma-1}}, \quad {\it as} \quad N\to\infty.$$
\end{thm}
{\bf Remark 1} This result coincides with the result that we got  in  the case of the periodic walk. There of course  
$$\gamma=\frac{1}{2L}\sum_{j=0}^{L-1} \frac{1}{p_j}.$$

\noindent
{\bf Remark 2} Our condition (i) is identical with Heyde's condition in (\ref{h}).

The organization  of the paper is the  following. In Section 2 we collect some  earlier results, which we need, and prove some lemmas. Section 3 contains the proof of Theorem 1.1. Section 4 contains the implications of Theorem 1.1, namely the local time results.
\noindent

\section{Preliminaries.}
 \renewcommand{\thesection}{\arabic{section}} \setcounter{equation}{0}
 \setcounter{thm}{0} \setcounter{lemma}{0}

First we are to redefine our random walk $\{{\mathbf C}(N);\,
N=0,1,2,\ldots\}$. It will be seen that the process described right below
is equivalent to that given in the Introduction.

To begin with, on a suitable probability space  consider two independent
simple symmetric (one-dimensional) random walks $S_1(\cdot)$, and
$S_2(\cdot)$. We may assume that on the same probability space we have a
double array of independent geometric random variables $\{G_i(j),\,
i\geq 1, j\in \mathbb{Z}\}$ which are independent from $S_1(\cdot)$, and
$S_2(\cdot),$ where $G_i(j)$ has the following geometric distribution
\begg
\mathbf{P}(G_i(j)=k)=2p_j(1-2p_j)^k=\alpha_j(1-\alpha_j)^k \,\, k=0,1,2,\ldots \quad j=0,\pm1,\pm2,\ldots \label{geo}
\endd
where  we used the notation $\alpha_j=2p_j,$ which will be convenient later.
We now construct our walk $\mathbf{C}(N)$ as follows. We will take all
the horizontal steps consecutively from $S_1(\cdot)$ and all the vertical
steps consecutively from $S_2(\cdot).$ First we will take some horizontal
steps from $S_1(\cdot)$, then exactly one vertical step from $S_2(\cdot),$
then again some horizontal steps from $S_1(\cdot)$ and exactly one vertical
step from $S_2(\cdot),$ and so on. Now we explain how to get the number of
horizontal steps on each occasion.
Consider our walk starting from the origin proceeding first horizontally
$G_1(0)$ steps (note that $G_1(0)=0$ is possible with probability
$2p_0$), after which it takes exactly one vertical step, arriving
either to the level $1$ or $-1$, where it takes $G_1(1)$ or
$G_1(-1)$ horizontal steps (which might be no steps at all) before
proceeding with another vertical step. If this step carries the walk to
the level $j$, then it will take $G_1(j)$ horizontal steps, if this is
the first visit to level $j,$ otherwise it takes $G_2(j)$ horizontal
steps. In general, if we finished the $k$ -th vertical step and arrived to
the level $j$ for the $i$-th time, then it will take $G_i(j)$
horizontal steps.

Let now $H_N,\, V_N$ be the number of horizontal and
vertical steps, respectively from the first $N$ steps of the
just described process. Formally 
\begin{eqnarray}
H_N&=&\#\{ k: \quad 1 \leq k\leq N, \,\,\, C_1(k)\neq C_1(k-1)\}\\
V_N&=& \#\{ k: \quad 1 \leq  k\leq N, \,\,\, C_2(k)\neq C_2(k-1)\}\\
\nonumber
\end{eqnarray}

Consequently, $H_N+V_N=N$, and
$$
\left\{{ \mathbf C}(N);\, N=0,1,2,\ldots\right\}=
\left\{(C_1(N),C_2(N));\, N=0,1,2,\ldots\right\}
$$
\begg
\stackrel{{d}}{=}\left\{(S_1(H_N),S_2(V_N));\, N=0,1,2,\ldots\right\},
\label{equ}
\endd
where $\stackrel{{d}}{=}$ stands for equality in distribution.

\smallskip

We will use   the following  theorem of Stenlund which will  be an important ingredient in our proof.

\bigskip
\noindent

{\bf Theorem  C }  \cite{ST} {\it Let $( Y_n)_{n\geq 0}$ be a Markov chain starting at $0$ with the following transition probabilities 

$${ \mathbf P}(Y_{n+1}=j+1 |Y_{n}=j)={ \mathbf P}(Y_{n+1}=j -1|Y_{n}=j)=p_j  \quad{\rm and} \quad   { \mathbf P}(Y_{n+1}=j |Y_{n}=j)=1-2p_j  $$
with some $0 <p_j\leq1/2$ for  all $j\in \mathbb{Z}. $ 
Suppose that 
\smallskip

\noindent
{\rm A1} $\quad$ The central limit theorem holds for $\frac{Y_n}{\sqrt{n}}$ for some $\sigma >0 $ :

 $$\lim_{n\to \infty}{ \mathbf P}(\frac{Y_n}{\sqrt{n}}\leq x)=\frac{1}{\sqrt{2\pi}\sigma}
\int_{-\infty}^x e^{-\frac{u^2}{2\sigma^2}}\,du=\int_{-\infty}^x \phi_{\sigma^2}(u)\,du.$$
\smallskip

\noindent
{\rm A2} $\quad$ There exists a $\mu>0$ for which 

\smallskip

$$\lim_{T\to \infty} \frac{1}{y-x} \int_x^y \frac{1}{p_{[Tu]}} du=\mu,  \qquad  {\it for\,\, all} \,\,(x,y)\in \mathbb R^2 $$
where $[u]$ means the largest integer $\leq u.$

\smallskip

\noindent
{\rm A3} \quad For the $n$-step transition probability  $(P^n)_{j,k}={ \mathbf P}(Y_n=k | Y_0=j) $ we have
$$\sup_{n\geq 0} \sup_{k \in \mathbb{Z}}\sqrt{n} ( P^n)_{0,k}<\infty.$$

\smallskip

\noindent
{\rm A4} \quad Suppose that we have

 $$ \sup_{n\geq 0}\sqrt{n} \sup_{N\geq n}\, \sum_{m=n}^N \left({ \mathbf P}(Y_{2m+2}=0 )-{ \mathbf P}(Y_{2m+1}=0 )\right)< \infty.$$

Then 
$$\lim_{n \to \infty} \sup_{x\in I}\left| p_{[\sqrt{n}x] }\sqrt{n}\,{ \mathbf P}(Y_n=[\sqrt{n}x])-\frac{1}{\mu} \phi_{\sigma^2}(x)\right |=0$$
on any compact set $I\in R.$ }

We will also need Lemma 9 from the above-mentioned paper 

\smallskip

\noindent
{\bf Lemma D}  Stenlund \cite{ST} {\it  If \, $\inf_{k\in {\mathbb Z}} \,p_k>0,$ there exists a constant  ${K}$ such  that

$$ (P^n)_{i,j}\leq \frac{K}{\sqrt{n}},\quad  (i,j) \in {\mathbb Z^2}, \quad n\geq 1.$$
 In particular, under  the condition of Lemma  {\rm D, A3} holds.}

\smallskip

We quote  the following remarks of Stenlund about his condition A4 in his Theorem C:

\smallskip

\noindent
{\bf Remark 3} If the walk is lazy, i.e. all  $p_j \leq 1/4,$ then  A4 is trivially satisfied. 
Simulation suggests that generally under the other conditions of Theorem C, A4 holds if $p_k \neq 1/2$ for at 
least one $k \in \mathbb {Z}.$ 
So A4 could well turn out  to be equivalent with the Markov chain being aperiodic.

\bigskip
Here we list a few facts about the simple symmetric walk and its local time. We define the local time of $S_n$ as
$$\xi(i,n)=\{\# j: 1\leq j \leq n: S_j=i \}.$$

\smallskip

\noindent
{\bf Lemma F} (F\"oldes and  R\'ev\'esz \cite{FR})
{\it Let $S_n=\sum_{i=1}^n  X_i$ be a simple symmetric walk, where 

${ \mathbf P}(X_i=1)={ \mathbf P}(X_i=-1)=1/2, \, i=1,2, \ldots$ and let $\xi(i,n)$ its local time. Let 

\[ I(i,\ell)=\left \{ \begin{array}{ll}
  1&\mbox {\,if \, $S_{\ell}=i \,or \, i+1$}\\
0 &\mbox {otherwise}\\
 \end{array}
 \right.\]

\noindent
and  $$D(i,n)=\xi(i+1,n)-\xi(i,n) -\sum_{\ell=0}^n I(i,\ell)X_{\ell+1}.$$
Then }
$$|D(i,n)|\leq 2.$$
Being $I(i,\ell)$ and $X_{\ell+1}$ independent and $\mathbf{E}(X_{\ell+1})=0$ we have

\smallskip
\noindent
{\bf Consequence}: 
 $$|\mathbf{E} (\xi(i+1,n)-\xi(i,n) )|\leq 2.$$

\bigskip

We need the following   well known result about the $n$-steps return time 

\smallskip
\noindent
{\bf Lemma G}  (see  e.g. \cite{RE}, Theorem 2.8 page  19)

$${ \mathbf P}(S_{2N}=0)={2N \choose N} \frac{1}{2^{2N}}=\frac{1}{\sqrt{\pi N}} 
\left(1+O\left(\frac{1}{N} \right) \right), \quad  {\it as} \quad N \to \infty.$$

\smallskip

For the second coordinate of our anisotropic walk we have

\begin{lemma} Under the conditions of Theorem {\rm1.1}
$$\lim _{N\to \infty}  \sqrt{N}{ \mathbf P}(C_2(2N)=0)=\frac{1}{4p_0\sqrt{\pi \gamma}}.$$
\end{lemma}
{\bf Proof:}
Observe that under the conditions of Theorem 1.1,\, $C_2(n)$ satisfies the conditions of 

\noindent
Theorem C. 
To see this, enough to note that Theorem A implies  A1  with $\sigma=\frac{1}{\sqrt{\gamma }}.$
A2 is a trivial consequence of (i) of Theorem 1.1 with $\mu=2\gamma.$ We have A3  
from Lemma D and condition (ii) of Theorem 1.1. Finally (iii) is identical to A4. Then we get by Theorem C for $x=0,$ that

$$\lim_{N\to \infty} p_0 \sqrt{2N}  { \mathbf P}(C_2(2N)=0)=\frac{1}{2\gamma} \phi_{\sigma^2}(0)=\frac{1}{2\sqrt{\gamma}\sqrt{2\pi}},$$ 
implying our lemma. 
$\Box$

We will need the following trivial lemma
\begin{lemma} Let  $ { \mathbf P}(G=k)=\alpha(1-\alpha)^k, \quad k=0,1,\ldots$ be a geometric random variable with success probability $\alpha.$ Let   $G^{(L)} $  be a truncated  geometric random variable, truncated  at some integer $L$  with distribution
$ { \mathbf P}(G^{(L)} =k)=\alpha(1-\alpha)^k, \quad k=0,1,\ldots L-1 $ and 
$ { \mathbf P}(G^{(L)} =L)=(1-\alpha)^L$. Then $Var (G^{(L)} )\leq\frac{2}{\alpha^2}.$
 
\end{lemma}
\noindent
{\bf Proof}: Clearly for any $L$, 

$$ \mathbf{E}(G^{(L)} )\leq \mathbf{E}(G) \qquad {\rm and} \qquad \mathbf{E}((G^{(L)} )^2)\leq \mathbf{E}(G^2).$$
So 
$$\mathbf{Var}(G^{(L)} )\leq \mathbf{E}((G^{(L)} )^2)\leq \mathbf{E}(G^2)= \mathbf{Var}(G)+(\mathbf{E}(G))^2=\frac{1-\alpha}{\alpha^2}+\frac{(1-\alpha)^2}{\alpha^2}\leq \frac{2}{\alpha^2}.$$
$\Box$

Define
\begg H_N^+= \sum_j\sum_{i=1}^{\xi_2(j,V_N)}G_i(j),  \label{hnp}
\endd
where $\{G_i(j) \quad i=1,2 \dots \quad j=0,\pm1,\pm2 ,\dots \}$ are independent random variables with distribution  (\ref{geo}).

 Then ${\mathbf E}(G_1(j))=\frac{1-\alpha_j}{\alpha_j}\,$ and 
 ${\mathbf Var}(G_1(j))=\frac{1-\alpha_j}{\alpha^2_j}.$

 Observe that $$ \sum_{i=1}^{\xi_2(j,V_N)}G_i(j)$$
 is the number of horizontal steps on the level $j$ in the first $N$ steps.  This is correct unless $j$ is the level, where the $N$-th step occurs. In this case, the last geometric random variable might be truncated. Thus the only difference between $H_N$ and $H_N^+$ is,  that the latter contains this last whole geometric random variable,  while $H_N$ might  contain only the possibly truncated  version of it.
 We will suppose throughout this paper, as in condition (ii) of Theorem 1.1, that $p_j\geq \omega$ for each $j$ with some small $\omega>0.$

 \begin{lemma} Under the conditions of Theorem  {\rm1.1}, for any $\delta>0$
\begg
{ \mathbf P}(H_N^+ -H_N>N^{\delta})\leq\exp(-cN^{\delta}).   \label{differ} 
\endd
\end{lemma}
{\bf Proof}: Clearly
$$
H_N^+ -H_N\leq\max_jG_{\xi_2(j,V_N)}(j).
$$
In the sequel we use the following notations
$$\sum_j=\sum_{\min_{0\leq k\leq V_N}S_2(k)\leq j\leq
\max_{0\leq k\leq V_N}S_2(k)}$$
 and
$$ \max_j=\max_{\min_{0\leq k\leq V_N}S_2(k)\leq j\leq
\max_{0\leq k\leq V_N}S_2(k)}
$$

\begin{eqnarray}
{ \mathbf P}(\max_j G_{\xi_2(j,V_N)}(j)>N^{\delta})&\leq& 
\sum_j { \mathbf P}(G_1(j)>N^{\delta})
\leq \sum_j (1-\alpha_j)^{N^{\delta}}
\\  \nonumber
&\leq&\sum_j (1-2\omega)^{N^{\delta}}
\leq(2N+1) \exp(-c_1 N^{\delta})
\leq \exp (-c N^{\delta})
  \label{bound}
\end{eqnarray}
with $c_1=-\ln (1-2\omega)$, and with some  $0<c<c_1.$
$\Box$

Introduce

\begg \frac{1}{j}\sum_{k=1}^j\frac{1}{p_k}=\kappa_j \qquad  \frac{1}{j}\sum_{k=1}^j\frac{1}{p_{-k}}=\beta_j .  \label{kappa} \endd
\begin{lemma} If $\ep_N= 1/N^{\rho}$
with  $0<\rho<1/4$, then 
$${ \mathbf P}\left(|H_N-\frac{\gamma-1}{\gamma}N|>N\ep_N\right)\leq { \mathbf P}\left(|H_N- \mathbf{E}(H_N)|>\frac{N\ep_N}{2}\right).$$

\end{lemma}

\noindent
{\bf Proof}:  We have from (\ref{hnp}) 

\begin{eqnarray}   \mathbf{E}(H_N^+)&=&
\sum_j 
 \mathbf{E}{\xi_2(j,V_N)} \frac{1-2p_j}{2p_j}= \sum_j 
 \mathbf{E}{\xi_2(j,V_N)} \frac{1}{2p_j}- \mathbf{E}(V_N). \label{ehaenp}
\end{eqnarray}
\label{fontos}
The following lines are  coming from Heyde \cite{H}

$$\sum_j \xi_2(j,V_N)\frac{1}{p_j}=\sum_{j=1}^\infty
\xi_2(j,V_N)(j\kappa_j-(j-1)\kappa_{j-1})+\sum_{j=1}^\infty
\xi_2(-j,V_N)(j\beta_j-(j-1)\beta_{j-1})+\xi_2(0,V_N)\frac{1}{p_0}$$
$$=\sum_{j=1}^\infty j \kappa_j(\xi_2(j,V_N)- \xi_2(j+1,V_N))+
\sum_{j=1}^\infty j \beta_j(\xi_2(-j,V_N)- \xi_2(-j-1,V_N))+
\xi_2(0,V_N)\frac{1}{p_0}$$
$$=\sum_{j=1}^\infty j(\kappa_j-2 \gamma)(\xi_2(j,V_N)-
\xi_2(j+1,V_N))+2\gamma \sum_{j=1}^\infty j(\xi_2(j,V_N)- \xi_2(j+1,V_N))$$
$$+\sum_{j=1}^\infty j(\beta_j-2\gamma)(\xi_2(-j,V_N)-
\xi_2(-j-1,V_N))+2\gamma \sum_{j=1}^\infty j(\xi_2(-j,V_N)-
\xi_2(-j-1,V_N))+\xi_2(0,V_N)\frac{1}{p_0}$$
$$=2\gamma \sum_{j=-\infty}^{\infty} \xi_2(j,V_N)+
\sum_{j=1}^{\infty} j(\kappa_j-2 \gamma)(\xi_2(j,V_N)
-\xi_2(j+1,V_N))$$
\begg +\sum_{j=1}^{\infty}
j(\beta_j-2\gamma)(\xi_2(-j,V_N)- \xi_2(-j-1,V_N))
+\xi_2(0,V_N)\left(\frac{1}{p_0}-2\gamma\right). \label{four}\endd

Using the notation of (\ref{kappa}),  we have under condition (i) of  Theorem 1.1  \begg 
j(|\kappa_j-2 \gamma|)\leq C j^{1-\eta}  \quad  {\rm and}\quad j(|\beta_j-2\gamma|)\leq C j^{1-\eta},
\quad j=1,2,\dots \label{Ce}
\endd
for some $C>0$ big enough.
Observe that
$$2\gamma \sum_{j=-\infty}^{\infty} \xi_2(j,V_N)=2\gamma V_N.$$
Then from Lemma F and (\ref{Ce}) we get as $\eta >1/2$ that as $N\to \infty$
$$\left| \mathbf{E}\left(\sum_{j=1}^{\infty} j(\kappa_j-2 \gamma)(\xi_2(j,V_N)
-\xi_2(j+1,V_N))\right)\right |\leq$$
$$ \mathbf{E}\left(\sum_{j=1}^{\infty} j|(\kappa_j-2 \gamma)||(\xi_2(j,V_N)
-\xi_2(j+1,V_N))|\right)= \mathbf{E}\left(\sum_{j=1}^{[N^{1/2+\delta}]}j|(\kappa_j-2 \gamma)||(\xi_2(j,V_N)
-\xi_2(j+1,V_N))|\right) $$
$$ +\mathbf{E}\left(\sum_{j=[N^{1/2+\delta}]+1}^Nj|(\kappa_j-2 \gamma)||(\xi_2(j,V_N)
-\xi_2(j+1,V_N))|\right)\leq 2\sum_{j=1}^{[N^{1/2+\delta}]}C j^{1-\eta}$$ 
$$+N\mathbf{P}(\max_{n\leq N} S_2(n)\geq N^{1/2+\delta}) \sum_{[N^{1/2+\delta}]}^N C j^{1-\eta}
=O(N^{(1/2+\delta)(2-\eta)})+\exp(-N^{c\delta})O(N^{3-\eta}),$$

\noindent
where $c > 0  $ and $\delta>0$ above are  small constants the values of which are unimportant.  In the last line of the above formula, we used the the facts that (see e.g. \cite{RE} page 21 formula (2.16) ) for any $\ep>0$ and $x_N=o(N^{1/6})$ 

$$\mathbf{P}\left(\frac{\max_{n\leq N}S_2(n)}{\sqrt{N}}\geq x_N\right)\leq (1+\ep )\frac{4}{\sqrt{2\pi}x_N}e^{-
x^2_N/2} \quad \quad  x_N>0. $$
Being $\eta>1/2$ we can select $0<\delta<1/6$ to be small enough that $(1/2+\delta)(2-\eta)\leq 3/4.$
 Thus we have
 
$$ O(N^{(1/2+\delta)(2-\eta)})+ \exp(-N^{c\delta})O(N^{4-\eta})=O(N^{3/4}).$$
 Similarly, the third sum in (\ref{four}) has the same expected value.
Finally for the last term (\ref{four}) we have 
 $$ \mathbf{E}\left(\xi_2(0,V_N)\left(\frac{1}{p_0}-2\gamma\right)\right)=O(N^{1/2+\delta}).$$
So we conclude that

$$ \mathbf{E}\left(\sum_{j}\frac{ \xi_2(j,V_N)}{p_j}\right)= 2\gamma  \mathbf{E}(V_N)+O(N^{3/4}).$$

\noindent
Then from (\ref{ehaenp})
$$ \mathbf{E}(H_N^+)=(\gamma-1) \mathbf{E}(V_N)+O(N^{3/4}).$$

\noindent
Observe that from  Lemma 2.3
$$ \mathbf{E}(H_N^+ -H_N)\leq N^{\delta}+N\exp(-cN^{\delta}).$$
Consequently
$$  \mathbf{E}(H_N)= \mathbf{E}(H_N^+)+O(N^{\delta}) $$  as well, implying that
$$E(H_N)=(\gamma-1)E(V_N)+O(N^{3/4}).$$
So 
$$N=\gamma  \mathbf{E}(V_N)+O(N^{3/4}) $$
 and 
 $$  \mathbf{E}(V_N)=\frac{1}{\gamma}N +O(N^{3/4}) \quad {\rm and} \quad  \mathbf{E}(H_N)=\frac{\gamma-1}{\gamma}N+O(N^{3/4}) .$$

\noindent
Thus $ \mathbf{E}(H_N-\frac{\gamma-1}{\gamma}N)=O(N^{3/4})$  implying that
$$ { \mathbf P}\left(|H_N-\frac{\gamma-1}{\gamma}N|>N\ep_N\right)\leq { \mathbf P}\left(|H_N- \mathbf{E}(H_N)|>\frac{N\ep_N}{2}\right)$$
for $\ep_N=1/ {N^{\rho}}$  with $\rho<1/4. $
$\Box$

\begin{lemma} Under the condition $p_j\geq \omega$ for $j=0,\pm1,\pm2 \ldots,$ and $\ep_N=1/ {N^{\rho}}$  with $0<\rho<1/4$ we have
$${ \mathbf P}\left(|H_N-\frac{\gamma-1}{\gamma}N|>N\ep_N\right)
 \leq \frac{2}{N\ep_N^2 \omega^2}.$$
  \end{lemma}

\noindent
{\bf Proof}: Observe from (\ref{geo}) that  $$\mathbf{Var}(G_i(j))=\frac{1-\alpha_j }{\alpha^2_j} \leq\frac{1}{4\omega^2}, \quad i=1,2 \dots \,\,,\,\,j=\pm1,\pm2 \ldots $$
Being $H_N^+$   the sum of $H_N$ and  an independent  truncated geometric random variable as in Lemma 2.2  and using that $V_N\leq N,$ we see that
$$\mathbf{Var}(H_N)\leq( N-1)\frac{1}{4\omega^2}+\frac{1}{2\omega^2}\leq \frac{N}{2\omega^2}.$$

Consequently, by Chebyshev inequality and Lemma 2.4, we get that 

$${ \mathbf P}\left(|H_N-\frac{\gamma-1}{\gamma}N|>N\ep_N\right)\leq { \mathbf P}\left(|H_N- \mathbf{E}(H_N)|>\frac{N\ep_N}{2}\right)
 \leq \frac{2}{N\ep_N^2 \omega^2}.$$
$\Box$

\section{Proof}
 \renewcommand{\thesection}{\arabic{section}} \setcounter{equation}{0}
 \setcounter{thm}{0} \setcounter{lemma}{0}

{\bf Proof of the Theorem 1.1.}  
Introduce the notation 
$$\gamma^*=\frac{\gamma-1}{\gamma}$$
\begin{eqnarray}
{ \mathbf P}(&{\mathbf C}&(2N)=(0,0))=
  \sum_{r=0} ^N { \mathbf P}\left ({\mathbf C}(2N)=(0,0)| H(2N)=2r){ \mathbf P}(H(2N)=2r\right) \nonumber \\
&=&\sum_{r=0} ^N { \mathbf P}(   {\mathbf C}(2N)=(0,0), S_1(2r)=0| H(2N)=2r){ \mathbf P}(H(2N)=2r) \nonumber\\
&=&\sum_{|r - \gamma^* N|<N\ep_N} { \mathbf P}({\mathbf C}(2N) =(0,0), S_1(2r)=0| H(2N)=2r){ \mathbf P}(H(2N)=2r)\nonumber\\
&+&\sum_{|r - \gamma^* N |\geq N\ep_N} { \mathbf P}({\mathbf C}(2N)=(0,0), S_1(2r)=0| H(2N)=2r){ \mathbf P}(H(2N)=2r)=I +II,
\end{eqnarray}
where $\ep_N>0$ is some small positive number.
We consider first the term I, which can be rewritten as 
\begin{eqnarray}
I= \sum_{|r - \gamma^* N|<N\ep_N}{ \mathbf P}(S_1(2r)=0)  { \mathbf P}(C_2(2N)=0  | H(2N)=2r){ \mathbf P}(H(2N)=2r).\nonumber
\end{eqnarray}
For  all $r$ satisfying  $|r - \gamma^* N|<N\ep_N,$ we use Lemma G and the inequalities: For $ \,
 0<x<1/2$
$$    \frac{1}{\sqrt{1-x}}\leq 1+x \quad {\rm and}\quad \frac{1}{\sqrt{1+x}}\geq1-x ,$$
to get
\begin{eqnarray}
{ \mathbf P}(S_1(2r)=0)&=&\frac{1}{\sqrt{\pi r}}\left (1+O\left(\frac{1}{r}\right)\right)\leq \nonumber
\frac{1}{\sqrt{\pi (\gamma^* N-N\ep_N)}}\left(1+O\left(\frac{1}{\gamma^* N-N\ep_N}\right)\right) \\ \label{remain1}
&=&\frac{1}{\sqrt{\pi \gamma^* N(1-\frac{\ep_N}{\gamma^*} )}}\left (1+O\left(\frac{1}{N}\right)\right)\leq
\frac{1}{\sqrt{\pi \gamma^* N}}\left (1+\frac{\ep_N}{\gamma^*}\right).
\end{eqnarray}
\noindent
Similarly,  we get that
\begin{eqnarray}
{ \mathbf P}(S_1(2r)=0)&=&\frac{1}{\sqrt{\pi r}} \left (1+O\left(\frac{1}{r}\right)\right)\geq 
\frac{1}{\sqrt{\pi( \gamma^* N+N\ep_N)}}\left(1+O\left(\frac{1}{\gamma^* N+N\ep_N}\right)\right) \nonumber\\
&=&\frac{1}{\sqrt{\pi \gamma^* N(1+\frac{\ep_N}{\gamma^*} )}}\left (1+O\left(\frac{1}{N}\right)\right)\geq
\frac{1}{\sqrt{\pi \gamma^* N}}\left (1-\frac{\ep_N}{\gamma^*}\right).\label{remaintwo}
\end{eqnarray}

Now let  $\ep_N=1/N^{\rho}$ with  some $0<\rho<1/4$ as in Lemma  2.5, then
\begin{eqnarray}\frac{1}{\sqrt{\pi \gamma^* N}}\left (1-\frac{\ep_N}{\gamma^*}\right)
\sum_{|r - \gamma^* N|<N\ep_N} { \mathbf P}(C_2(2N)=0,  | H(2N)=2r){ \mathbf P}(H(2N)=2r)\leq I \nonumber\\
\leq \frac{1}{\sqrt{\pi \gamma^* N}}\left (1+\frac{\ep_N}{\gamma^*}\right) \sum_{|r - \gamma^* N|<N\ep_N} { \mathbf P}(C_2(2N)=0,  | H(2N)=2r){ \mathbf P}(H(2N)=2r). \label{extend}
\end{eqnarray}

Now observe that
 \begin{eqnarray} \sum_{|r - \gamma^* N|\geq  N\ep_N}&{ \mathbf P}&( C_2(2N)=0| H(2N)=2r){ \mathbf P}(H(2N)=2r) 
\nonumber\\
&=&\sum_{|r - \gamma^* N|\geq N\ep_N}{ \mathbf P}( C_2(2N)=0, H(2N)=2r)\leq \sum_{|r - \gamma^* N|\geq N\ep_N}{ \mathbf P}(  H(2N)=2r) \nonumber\\
&=&{\bf P}(  |H(2N)-2\gamma^*N | \geq 2N\ep_{2N})\leq \frac{1}{N\ep^2_{2N}\omega^2}\leq \frac{c}{N^{1-2\rho}}, 
\label{firstrem}
\end{eqnarray}
where we used Lemma 2.5 in the last line. Here $c$ is an unimportant constant depending  only on $\omega .$ 
Extending  now the summation  for all $r$ in (\ref{extend}) and observing that $$\sum_r { \mathbf P}( C_2(2N)=0| H(2N)=2r){ \mathbf P}(H(2N)=2r)={ \mathbf P}(C_2(2N)=0),$$
we  conclude by using (\ref {extend}) and (\ref{firstrem})  that

 \begin{eqnarray} \frac{1}{\sqrt{\pi \gamma^* N}}\left (1-\frac{\ep_N}{\gamma^*}\right)
\left(  { \mathbf P}( C_2(2N)=0) -\frac{c}{N^{1-2\rho}}\right)\leq I\\ \nonumber
\leq \frac{1}{\sqrt{\pi \gamma^* N}}\left (1+\frac{\ep_N}{\gamma^*}\right) \left(  { \mathbf P}( C_2(2N)=0) +\frac{c}{N^{1-2\rho}}\right).
\end{eqnarray}
Thus we can get that
\begin{eqnarray} 
I=\frac{\sqrt \gamma}{\sqrt{\pi (\gamma-1)N}}
{ \mathbf P}( C_2(2N)=0)\left(1+O(\ep_N)\right)+O\left(\frac{1}{N^{3/2-2\rho}}\right). \label{c2} 
\end{eqnarray}

Then using (\ref{c2}) and Lemma 2.1, we get that
\begg I =\frac{1}{4Np_0\pi\sqrt{\gamma-1}}\left(1+O(\ep_N)\right)+O\left(\frac{1}{N^{3/2-2\rho}}\right).  \label{ime I} 
\endd

Then what remains to consider is
\begin{eqnarray*} II= \sum_{|r - \gamma^* N| \geq N\ep_N}{ \mathbf P}(C(2N)=(0,0), S_1(2r)=0| H(2N)=2r){ \mathbf P}(H(2N)=2r)\\
=\sum_{|r - \gamma^* N| \geq N\ep_N}{ \mathbf P}( S_1(2r)=0, S_2(2N-2r)=0| H(2N)=2r){ \mathbf P}(H(2N)=2r).
\end{eqnarray*}
\noindent
If
$$|r - \gamma^* N| \geq N\ep_N,$$
then either $r\geq N(\gamma^*+\ep_N)$ or   $r\leq N(\gamma^*-\ep_N).$ In this latter case $N-r\geq N(1-\gamma^*+\ep_N))$
(observe that $1-\gamma^*=1/\gamma$ is always positive).
This implies that in the first case ${ \mathbf P}( S_1(2r)=0)= O \left(\frac{1}{\sqrt{N}}\right)$    
and in the second case ${ \mathbf P}( S_2(2N-2r)=0)=O\left(\frac{1}{\sqrt{N}}\right).$
Thus

\begin{eqnarray}& II&= \sum_{|r - \gamma^* N| \geq N\ep_N} { \mathbf P}( S_1(2r)=0, S_2(2N-2r)=0| H(2N)=2r){ \mathbf P}(H(2N)=2r) \nonumber\\
&\leq&O \left(\frac{1}{\sqrt{N}}\right)\left(\sum_{r\geq N(\gamma*+\ep_N)} { \mathbf P}(  S_2(2N-2r)=0, H(2N)=2r)
+\sum_{r \leq N(\gamma*-\ep_N)} { \mathbf P}( S_1(2r)=0,  H(2N)=2r)\right)\nonumber\\
&\leq&O \left(\frac{1}{\sqrt{N}}\right)\left(\sum_{r\geq N(\gamma*+\ep_N)} { \mathbf P}(  H(2N)=2r)
+\sum_{r \leq N(\gamma*-\ep_N)} { \mathbf P}(  H(2N)=2r)\right) \nonumber \\
   \leq &O& \left(\frac{1}{\sqrt{N}}\right) \sum_{|r - \gamma^* N| \geq N\ep_N}{ \mathbf P}(  H(2N)=2r)=O \left(\frac{1}{\sqrt{N}}\right)
{ \mathbf P}(|H(2N)-2\gamma^* N|>2N\ep_N) \nonumber\\
&=&O\left(\frac{1}{N^{3/2}\ep^2_N}\right)=O\left(\frac{1}{N^{3/2-2\rho}}\right) \label{ime II}
\end{eqnarray}
by Lemma 2.5,  whenever  $\ep_N=\frac{1}{N^\rho}$ with $\rho< 1/4.$
By (\ref{ime I}) and (\ref{ime II})   selecting $\rho<1/4$ we get that the second term in I  and  term II are negligible to the main term in I, which gives the asymptotic probability of the $2N$ -step return time to zero. $\Box$

\section{Local time}
 \renewcommand{\thesection}{\arabic{section}} \setcounter{equation}{0}
 \setcounter{thm}{0} \setcounter{lemma}{0}

The  immediate  consequence of our Theorem 1.1 
is, that under the conditions of Theorem 1.1.  the truncated Green function $g(\cdot)$ is given by
$$
g(N)=\sum_{k=0}^N \mathbf{P}(\mathbf{C}(k)=(0,0))\sim \frac{\log
N}{4p_0\pi\sqrt{\gamma-1}},\qquad N\to\infty,
$$
which implies that our anisotropic random walk in this case is recurrent
and also Harris recurrent.

Now define the local time by
$$
\Xi((k,j),N)=\sum_{r=1}^N I\{\mathbf{C}(r)=(k,j)\}.
$$ 
To get the next result we need the invariant measure, which is defined as 

$$
\mu(A)=\sum_{(k,j)}\mu(k,j)\mathbf{P}(\mathbf{C}(N+1)\in
A|\mathbf{C}(N)=(k,j)).
$$
In our case
$$
\mu(k,j)=\mu(k+1,j)\left(\frac12-p_j\right)
+\mu(k-1,j)\left(\frac12-p_j\right)
+\mu(k,j+1)p_{j+1}+\mu(k,j-1)p_{j-1}.
$$
It is easy to see that
$$
\mu(k,j)=\frac1{p_j}, \quad (k,j)\in\mathbb{Z}^2
$$
satisfies this equation, so it is an invariant measure.
In the case when the random walk is (Harris) recurrent, then we have (cf.
e.g. Chen \cite{CX})
$$
\lim_{n\to\infty}\frac{\Xi((k_1,j_1),n)}{\Xi((k_2,j_2),n)}=
\frac{\mu(k_1,j_1)}{\mu(k_2,j_2)}\quad {a.s.}
$$
 Hence
$$
\lim_{n\to\infty}\frac{\Xi((0,0),n)}{\Xi((k,j),n)}=\frac{p_j}{p_0}
\quad {a.s.}
$$
for fixed $(k,j)$.

It follows from Darling and Kac \cite{DK} that we have exponential
limiting distribution:
\begin{corollary}
$$
\lim_{n\to\infty}\mathbf{P}\left(\frac{\Xi((0,0),n)}{g(n)}\geq x\right)=
\lim_{n\to\infty}
\mathbf{P}\left(\frac{4p_0\pi\sqrt{\gamma-1}\, \Xi((0,0),n)}{\log n}\geq
x\right)=e^{-x},\quad x\geq 0.
$$
\end{corollary}

For limsup result we have (cf. Chen \cite{CX}):
\begin{corollary}
$$
\limsup_{n\to\infty}\frac{\Xi((0,0),n)}{\log n\log\log\log n}
=\frac{1}{4p_0\pi\sqrt{\gamma-1}}\quad{a.s.}
$$
\end{corollary}

For moderate and large deviations and functional limit laws for the local
time see Cs\'aki et al. \cite{CRR}, which was extended by Gantert and
Zeitouni \cite{GZ}. In our case, the functional limit theorem reads as
follows: Let ${\cal M}$ be the set of functions $m(x)$, $0\leq x\leq 1$
which are non-decreasing, right-continuous on $[0,1)$ and left-continuous
at $x=1$, equipped with weak topology, induced by L\'evy metric.
Furthermore, let ${\cal M}^*$ be the subset of ${\cal M}$ with $m(0)=0$
and
$$
\int_0^1\frac{dm(x)}{x}\leq 1.
$$
\begin{corollary}
Let $t(n,x)\in{\cal M}$ be a sequence of functions such 
$$
\lim_{n\to\infty}\frac{\log t(n,x)}{\log n}=x
$$
for all $0\leq x\leq 1$. Put
$$
f_n(x)=\frac{4p_0\pi\sqrt{\gamma-1}\, \Xi(0,t(n,x))}{\log n\log\log\log n}
$$
Then almost surely, the set of limit points of
$\{f_n(x),\, 0\leq x\leq 1\}_{n\geq 16}$ is ${\cal M}^*$.
\end{corollary}

\noindent
{\bf Remark:} For example  $t(n, x)=n^x$ satisfy the condition of the Corollary.

\section*{Data Availability}

Data sharing is not applicable to this article as no dataset were generated or analyzed during the current study.

\section*{Conflict of interest}

The  authors declare that they have no conflict of interest.

\end{document}